\documentclass[10pt]{amsart}

\usepackage{graphicx}

\usepackage[USenglish]{babel}
\usepackage{latexsym}
\usepackage{amsmath}
\usepackage{amsfonts}
\usepackage{amssymb}
\usepackage{esint}
\usepackage[all]{xy}
\renewcommand{\a }{\alpha }
\renewcommand{\b }{\beta }

\newcommand{\D }{\Delta }

\newcommand{\M}{\mathcal{M}}

\newcommand{\g }{\gamma}

\renewcommand{\L }{\Lambda }

\newcommand{\n }{\nabla }

\newcommand{\Sg }{\Sigma}

\newcommand{\ov}{\overline}

\newcommand{\wtilde }{\widetilde}

\newcommand{\be}{\begin{equation}}
\newcommand{\ee}{\end{equation}}

\newcommand{\R}{\mathbb{R}}

\newcommand{\Z}{\mathbb{Z}}

\newcommand{\N}{\mathbb{N}}
\newcommand{\no}{\noindent}

\newtheorem{theorem}{Theorem}[section]
\newtheorem{proposition}[theorem]{Proposition}
\newtheorem{example}[theorem]{Example}

\newcommand{\bpr}{\begin{proposition}}
\newcommand{\epr}{\end{proposition}}
\newcommand{\bex}{\begin{example}\rm}
\newcommand{\eex}{\end{example}}

\begin{document}

\newtheorem{lem}{Lemma}[section]
\newtheorem{pro}[lem]{Proposition}
\newtheorem{thm}[lem]{Theorem}
\newtheorem{rem}[lem]{Remark}
\newtheorem{cor}[lem]{Corollary}
\newtheorem{df}[lem]{Definition}

\title[The Mean Field equation on compact surfaces]{A note on a multiplicity result for the mean field equation on compact surfaces}

\author {Aleks Jevnikar}

\address{SISSA, via Bonomea 265, 34136 Trieste (Italy).}


\email{ajevnika@sissa.it}



\begin{abstract}
We are concerned with the following class of equations with exponential nonlinearities on a compact surface $\Sg$:
$$
  - \D u = \rho_1 \left( \frac{h \,e^{u}}{\int_\Sg
      h \,e^{u} \,dV_g} - \frac{1}{|\Sg|} \right) - \rho_2 \left( \frac{h \,e^{-u}}{\int_\Sg
      h \,e^{-u} \,dV_g} - \frac{1}{|\Sg|} \right),
$$
which describes the mean field equation of the equilibrium turbulence with arbitrarily signed vortices. Here $h$ is a smooth positive
function and $\rho_1, \rho_2$ two positive parameters. 

We provide the first multiplicity result for this class of equations by using Morse theory.
\end{abstract}

\maketitle

\section{Introduction}

\no In this paper we consider the following mean field equation
\begin{equation}
-\D u = \rho_1 \left(\frac{h\,e^u}{\int_\Sigma h\,e^u \,dV_g} - \frac{1}{|\Sigma|}\right) - \rho_2 \left(\frac{h\,e^{-u}}{\int_\Sigma h\,e^{-u} \,dV_g} - \frac{1}{|\Sigma|} \right) \hspace{0.3cm} \mbox{on $\Sigma$}, \label{eq}
\end{equation}
where $\D = \D_g$ is the Laplace-Beltrami operator, $\rho_1, \rho_2$ are two non-negative parameters, $h:\Sigma \rightarrow \mathbb{R}$ is a smooth positive function and $\Sigma$ is a compact orientable surface without boundary with Riemannian metric $g$ and total volume $|\Sigma|$. Throughout the paper, for the sake of simplicity we assume that $|\Sigma|=1$.

Equation \eqref{eq} arises in mathematical physics as the mean field equation of the equilibrium turbulence with arbitrarily signed vortices and it was first introduced by Joyce and Montgomery \cite{joy-mont} and  by Pointin and Lundgren \cite{point-lund}. These vortices are composed of positive and negative intensities with the same value, where $u$ and $\rho_1/\rho_2$ are associated with the stream function of the fluid and the ratio of the numbers of the signed vortices, respectively. Later, several authors worked on this model; we refer for example to \cite{cho,lio, mar-pul,new,oh-su} and the references therein. The case $\rho_1 = \rho_2$ has a close relationship with geometry and is related to the study of constant mean curvature surfaces, see \cite{wen1,wen2}.

\medskip

Problem \eqref{eq} is variational, and solutions correspond to critical points of the Euler-Lagrange functional $J_\rho:H^1(\Sg) \to \R$, $\rho = (\rho_1,\rho_2)$, given by
\begin{eqnarray} \label{func}
J_\rho(u) &=& \frac{1}{2}\int_\Sigma |\nabla_g u|^2 \,dV_g + \rho_1 \left( \int_\Sigma u \,dV_g - \log\int_\Sigma h(x)\,e^u \,dV_g \right) +   \nonumber\\
                   & & -\,\, \rho_2 \left( \int_\Sigma u \,dV_g + \log\int_\Sigma h(x)\,e^{-u} \,dV_g \right).
\end{eqnarray}
A fundamental tool in dealing with this kind of functionals is the Moser-Trudinger inequality \eqref{eq:mt} and its version for the two-parameter case obtained in \cite{oh-su}:
\begin{equation} \label{mt2}
8\pi \left( \log\int_\Sigma e^{u-\ov u} \,dV_g  + \log\int_\Sigma e^{-u+\ov u} \,dV_g \right) \leq \frac 12 \int_\Sigma |\nabla_g u|^2 \,dV_g + C_\Sigma, 
\end{equation}
where $\ov{u}$ denotes the average of $u$. It follows directly that the functional $J_\rho$ is bounded from below and coercive whenever both $\rho_1$ and $\rho_2$ are less than $8\pi$. Therefore, one deduces the existence of a solution by minimization technique. On the other hand, when on of the $\rho_i$ exceeds the value $8\pi$ the functional is unbounded from below and the problem gets more involved.

\medskip

To describe the features of the problem and the general strategy to attack this kind of equations, it is first convenient to discuss its one-parameter counterpart, namely the following standard Liouville equation:
\begin{equation} \label{eq2}
  - \D u = \rho \left( \frac{h \,e^{u}}{\int_\Sg
      h \,e^{u} \,dV_g} - 1 \right).
\end{equation} 
Equation \eqref{eq2} concerns the problem in conformal geometry of prescribing the Gauss curvature of a surface, see \cite{bah-cor, cha, cha2, li, scho-zha}. More precisely, letting $\tilde g= e^{2v}g$, the Laplace-Beltrami operator of the deformed metric is given by $\D_{\tilde g} = e^{-2v}\D_g$ and the change of the Gauss curvature is ruled by
$$
	-\D_g v = K_{\tilde g} e^{2v} - K_g,
$$  
where $K_g$ and $K_{\tilde g}$ are the Gauss curvatures of $(\Sg, g)$ and of $(\Sg, \tilde g)$ respectively. Another motivation for the study of \eqref{eq2} is in mathematical physics as it models the mean field equation of Euler flows, see \cite{cal, kies}. This problem has been studied widely by lots of authors and there are by now many results regarding existence, blow-up analysis, compactness of solutions, etc, see \cite{djlw,djadli,mal,tar}.

An important feature of problem \eqref{eq2} is the lack of compactness, as its solutions might blow-up. In this case a quantization phenomenon was observed  in \cite{bre-mer, li1,li-shaf}. Indeed, a blow-up point $\ov{x}$ for a sequence $(u_n)_n$ of solutions relatively to $(\rho_n)_n$, i.e. there exists a sequence $x_n \to \ov x$ such that $u_n(x_n) \to +\infty$ as $n \to +\infty$, satisfies 
\begin{equation} \label{quant}
\lim_{r \to 0} \lim_{n \to + \infty} \rho_n \frac{\int_{B_r(\ov{x})} h \, e^{u_n} \,dV_g}{\int_\Sg h \, e^{u_n} \,dV_g} = 8 \pi.
\end{equation}
Somehow, each blow-up point has a quantized local mass. Furthermore, the limit profile of solutions is close to a {\em bubble}, namely a function $U_{\lambda,\ov x}$ defined as 
$$
	U_{\lambda,\ov x}(y) = \log \left( \frac{4\lambda}{\bigr( 1 + \lambda \,d(\ov x,y)^2 \bigr)^2} \right),
$$
where $y\!\in\!\Sg,\, d(\ov x,y)$ stands for the geodesic distance and $\lambda$ is a large parameter. In other words, the limit function is the logarithm of the conformal factor of the stereographic projection from $S^2$ onto $\R^2$, composed with a dilation.

In the general case when $\rho_2\neq 0$, namely for problem \eqref{eq}, the refined blow-up analysis is not carried out in full depth. Still, one can show that equation \eqref{eq} inherits some character from the Liouville case. In fact, in \cite{jwyz} the authors proved an analogous quantization property; for a blow-up point $\ov x$ and a sequence $(u_n)_n$ of solutions relatively to $(\rho_{1,n}, \rho_{2,n})$ one gets
$$
    \lim_{r \to 0} \lim_{n \to + \infty} \rho_{1,n} \frac{ \int_{B_r(\ov x)} h \, e^{u_n} \, dV_g }{ \int_\Sg h\, e^{u_n}\, dV_g } \in 8 \pi \N, \quad
    \lim_{r \to 0} \lim_{n \to + \infty} \rho_{2,n} \frac{ \int_{B_r(\ov x)} h \, e^{-u_n} \, dV_g }{ \int_\Sg h\, e^{-u_n}\, dV_g } \in 8 \pi \N.
$$
Moreover, it is possible to show that the case of multiples of $8 \pi$ indeed occurs, see \cite{es-we, gr-pi}. We let now
\begin{equation} \label{lambda}
	\L = (8 \pi \N \times \R) \cup (\R \cup 8 \pi \N).
\end{equation}
Combining the local volume quantization with some further analysis, see \cite{bat-man}, we get that the set of solutions is compact for $\rho_i$ bounded away from multiples of $8\pi$. This is the main reason why one has to restrict himself to parameters $(\rho_1,\rho_2) \notin \L$. In fact, in order to run the variational methods some compactness property is required, usually the Palais-Smale condition. Unfortunately, it is not known whether the latter holds or not for this equation. However, there is a way around it using the monotonicity argument from \cite{struwe} jointly with the compactness result.

\medskip

We briefly illustrate here the role played by the study of sublevels of the energy functional in the existence issue. Let us consider first the Liouville case \eqref{eq2}, with associated functional
\begin{equation} \label{fun-sca}
	I_\rho(u) = \frac{1}{2}\int_\Sigma |\nabla_g u|^2 \,dV_g + \rho\left(\int_\Sigma u \,dV_g - \log \int_\Sigma h(x)\,e^u \,dV_g\right). 
\end{equation}
From the standard Moser-Trudinger inequality
\begin{equation}\label{eq:mt}
   8\pi  \log \int_\Sigma e^{u - \ov{u}} \,dV_g  \leq \frac 12 \int_\Sigma |\n u|^2 dV_g + C_{\Sigma,g}, 
\end{equation}
we get boundedness and coercivity provided $\rho<8\pi$. For larger values of the parameter, the first step was done in \cite{cl3}, where an improved Moser-Trudinger is presented; roughly speaking, the more the function $e^u$ is spread over the surface $\Sg$, the better is the constant in the inequality and, as a consequence, one gets new lower bounds on the functional \eqref{fun-sca}. Improving this result, in \cite{djadli} and \cite{mal} the authors were able to deduce a general existence result. Basically, if $\rho<8(k+1)\pi$, $k\in\N$, and if $I_\rho(u)$ is large negative, i.e. lower bounds fail, $e^u$ has to be concentrated around at most $k$ points of $\Sg$. To represent this scenario it is then natural to consider the family of unit measures $\Sg_k$ which are supported in at most $k$ points of $Σ$, known as \emph{formal barycenters} of $\Sg$ of order $k$:
\begin{equation}\label{sigk}
	\Sg_k = \left\{ \sum_{i=1}^k t_i\delta_{x_i} \, : \, \sum_{i=1}^k t_i=1,t_i\geq 0,x_i\in\Sg,\forall\,i=1,\dots,k \right\}.
\end{equation}
The authors indeed proved a homotopy equivalence between the latter set and the low sublevels of $I_\rho$. The existence of solutions follows then from the non contractibility of $\Sg_k$ and suitable min-max schemes.

\medskip

Concerning the general case \eqref{eq}, the semi-coercive case $\rho_1\in(8k\pi, 8(k+1)\pi)$, $k\in\N$ and $\rho_2<8\pi$ was considered in \cite{zhou}. The author exploited the condition $\rho_2<8\pi$ to characterize the low sublevels of $J_\rho$ by means of the component $e^u$ only, which has the same concentration
behavior which occurs in the one-parameter case \eqref{eq2}. 

For parameters above the threshold value $(8\pi, 8\pi)$ the existence problem gets more involved and has still to be examined in depth due to the non trivial interaction of the two components $e^u$ and $e^{-u}$. It turns out that there is some analogy between this problem and the Toda system of Liouville equations arising from Chern-Simons theory. 

In the spirit of \cite{mr}, the first step was done in \cite{jev}, where the author derived an existence result for the first non trivial interval, i.e. $(\rho_1, \rho_2)\in(8\pi, 16\pi)^2$. The proof relies on an improved Moser-Trudinger inequality: one can indeed show that when both $e^u$ and $e^{-u}$ concentrate around the same point and with the same \emph{rate}, the constant in the left-hand side of \eqref{mt2} can be basically doubled. 

The general case with $(\rho_1, \rho_2)\notin \L$ was then considered in \cite{todatori} under the assumption that the surface $\Sg$ is not homeomorphic to $S^2$. The strategy goes as follows; exploiting improved Moser-Trudinger inequalities it is possible to show that if $\rho_1 < 8(k + 1)\pi$ and $\rho_2 < 8(l + 1)\pi$, $k, l \in \N$, then either $e^{u}$ is close to $\Sg_k$ or $e^{-u}$ is close to $\Sg_l$ in the distributional sense. This alternative can be expressed by means of the \emph{topological join} of $\Sg_k$ and $\Sg_l$. We recall that, given two sets $A$ and $B$, the join $A*B$ is defined as the family of elements of the form
\begin{equation}\label{join}
 A*B = \frac{ \bigr\{ (a,b,s): \; a \in A,\; b \in B,\; s \in [0,1]  \bigr\}}R,
\end{equation}
where $R$ is an equivalence relation such that
$$\begin{array}{l}
  (a_1, b,1) \stackrel{R}{\sim} (a_2,b, 1)  \quad \forall a_1, a_2 \in A, b \in B; \\
  (a, b_1,0)  \stackrel{R}{\sim} (a, b_2,0) \quad \forall a \in A, b_1, b_2 \in B.
\end{array}$$
Roughly speaking, the join parameter $s$ expresses which of the above alternatives is more likely to be fulfilled.

To minimize the interaction of the two components $e^u$ and $e^{-u}$ the assumption on $\Sg$ is needed. One can indeed construct two disjoint simple non-contractible curves $\gamma_1, \gamma_2$ such that $\Sg$ retracts on each of them through continuous maps $\Pi_1, \Pi_2$. Taking into account the retractions $\Pi_i$, starting from $\Sg_k * \Sg_l$, one can restrict himself to targets in $(\gamma_1)_k*(\gamma_2)_l$ only. The final step is then to gain some non trivial topological informations of the low sublevels of $J_\rho$ in terms of $(\gamma_1)_k*(\gamma_2)_l$.

\medskip

There are also some recent results concerning a Leray-Schauder degree theory approach for both the mean field equation and the Toda system, see \cite{jev2, mrd}.

\

The goal of this paper is to present the first multiplicity result for this class of equations. 
\begin{thm}\label{teo}
Let $\rho_1 \in (8k\pi, 8(k + 1)\pi)$ and $\rho_2 \in (8l\pi, 8(l + 1)\pi)$, $k, l \in \N$ and let $\Sg$ be a compact surface with genus $g(\Sg)>0$. Then, for a generic choice of the metric $g$ and of the function $h$ it holds
$$
 \#\bigr\{ \mbox{solutions of \eqref{eq}} \bigr\} \geq \binom{k+g(\Sg)-1}{g(\Sg)-1} \binom{l+g(\Sg)-1}{g(\Sg)-1}. 
$$
\end{thm}
Here, by generic choice of $(g,h)$ we mean that it can be taken in an open dense subset of $\M^2 \times C^{\,2}(\Sg)^+$, where $\M^2$ stands for the space of Riemannian metrics on $\Sg$ equipped with the $C^{\,2}$ norm, see Proposition \ref{non-deg}.  

\medskip

The proof is carried out by means of the Morse theory in the same spirit of \cite{bar-mar-mal} and \cite{batt}, where the problem of prescribing conformal metrics
on surfaces with conical singularities and the Toda system are considered, respectively. The argument is based on the analysis developed in \cite{todatori}: in particular we will exploit the topological descriptions of the low sublevels of $J_\rho$ to get a lower bound on the number of solutions to \eqref{eq}. It will turn out indeed that the high sublevels of $J_\rho$ are contractible, while the low sublevels carry some non trivial topology. We will finally apply the weak Morse inequalities to deduce the estimate on the number of solutions by means of the latter change of topology. Somehow, one expects that the more the topology of the surface $\Sg$ is involved, the higher is the number of solutions. In fact, we will exploit the genus of $\Sg$ to describe the topology of low sublevels of $J_\rho$ by means of some bouquet of circles, see Lemma \ref{bouquet} and Proposition \ref{pro1}. In this way we will capture the topological informations of $\Sg$ and provide a better bound on the number of solutions to \eqref{eq}. 

\medskip

The paper is organized as follows: in Section 2 we introduce some notations and we collect some known results we will use later on. In particular, we focus first on a compactness resul of equation \eqref{eq} and we introduce a deformation lemma for the functional $J_\rho$ in order to use Morse arguments. The second part is concerned with a classic result in Morse theory: the Morse inequalities. In Section 3 we finally prove the main result of Theorem \ref{teo}.

\section{Preliminaries}

\no We give here some notation and known results which we will use throughout the paper.

\subsection{Notation} \label{not}

\

\medskip

\no The genus of the surface $\Sg$ will be denoted by $g(\Sg)$. The space of Riemannian metrics on
$\Sg$ equipped with the $C^2$ norm  will be indicated by $\mathcal M^2$. The symbol $B_r(p)$ stands for the open metric ball of radius $r$ and centre $p\in \Sg$. Given a function $u\in L^1(\Sg)$, the average of $u$ is defined by
$$
	\ov{u} = \frac{1}{|\Sg|}\int_\Sg u \,dV_g.
$$
The sublevels of the functional $J_\rho$ will be denoted by
$$
	J_\rho^a = \bigr\{ u\in H^1(\Sg) : J_\rho(u)\leq a \bigr\}.
$$
The sign $\simeq$ will refer to homotopy equivalence, while $\cong$ will stand for homeomorphisms between topological spaces or isomorphisms between groups. The identity map on a space $X$ will be indicated by Id$_{X}$. 

Given $q\in \N$ and a topological space $X$, we will denote by $H_q(X)$ its $q$-th homology group with
coefficient in $\Z$. For a subspace $A\subseteq X$ we write $H_q(X,A)$ for the $q$-th relative homology group of $(X,A)$. We will denote by $\wtilde H_q(X)$ the reduced $q$-th homology group, i.e. $H_0(X)=\wtilde H_0(X) \oplus \Z$ and $H_q(X)=\wtilde H_q(X)$ for all $q>0$.

The $q$-th Betti number of $X$ will be indicated by $\b_q(X)$, namely $\b_q(X)= \mbox{rank } (H_q(X))$, while $\wtilde\b_q(X)$ will correspond to the rank of the reduced homology group.

The letter $C$ will stand for large constants which are allowed to vary among different
formulas. To stress the dependence of the constants on some parameter, we add subscripts to $C$.
\

\subsection{Compactness result and a Deformation Lemma}

\

\medskip

\no We state now the compactness result of the set of solutions of equation \eqref{eq}. Recall the definition of the set $\L$ given in \eqref{lambda}. As mentioned in the Introduction, the blow-up phenomenon yields a quantization property of the local volume and with some standard analysis, see \cite{bat-man}, we deduce the following:
\begin{thm}(\cite{bat-man, jwyz}) \label{compt}
For $(\rho_1,\rho_2)$ in a fixed compact set of $\R^2\setminus \L$ the family of solutions to \eqref{eq} is uniformly bounded in $C^{2,\a}$ for some $\a>0$.
\end{thm}
We will need the latter compactness property to bypass the Palais-Smale condition, since it is not know whether it holds or not for this class of equations. More precisely, one can adapt the strategy in \cite{lucia}, where a deformation lemma for the Liouville equation \eqref{eq2} was presented, for our framework, see also \cite{mal1, struwe2}. One has the next alternative: either there exists a critical point of the functional $J_\rho$ inside some interval or there is a deformation retract between the relative sublevels. Recall the notation for the sublevels $J_\rho^a$ given in Subsection \ref{not}.  
\begin{lem}
If $\rho=(\rho_1,\rho_2)\notin \L$ and if $a<b\in\R$ are such that $J_\rho$ has no critical levels inside the interval $[a,b]$, then $J_\rho^a$ is a deformation retract of $J_\rho^b$.
\end{lem} 
Here, by deformation retract of a space $X$ onto some subspace $A\subseteq X$ we mean a continuous map $R:[0,1]\times X \to X$ such that $R(t,a)=a$ for all $(t,a)\in [0,1]\times A$ and such that the final target of $R$ is contained in $A$, i.e. $R(1,\cdot)\in A$.

\medskip

\begin{figure}[h]
\centering
\includegraphics[width=0.6\linewidth]{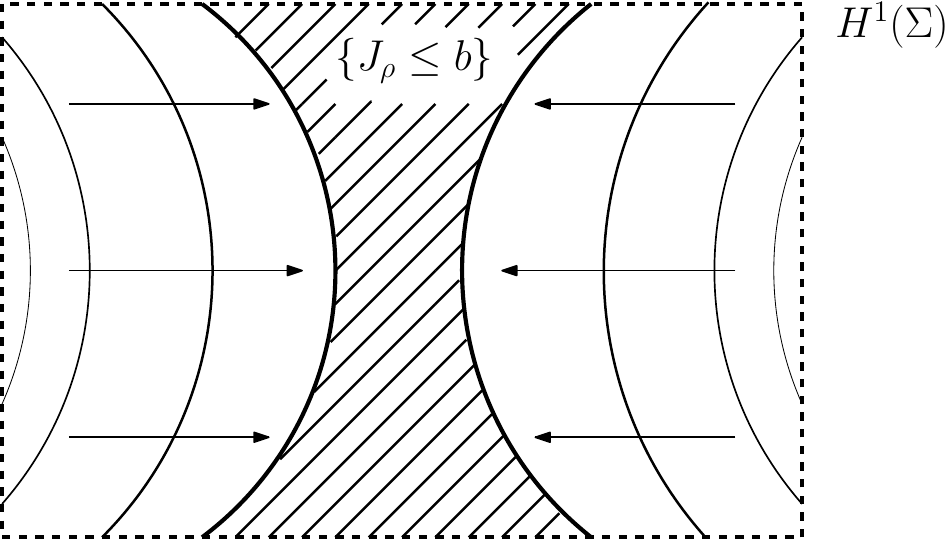}
\caption{}
\label{fig:alti}
\end{figure}

Notice now that by the compactness result of Theorem \ref{compt} it follows that $J_\rho$ has no critical points above some high level $b\gg 0$. Therefore, one can obtain a deformation retract of the whole Hilbert space $H^1(\Sg)$ onto the sublevel $J_\rho^b$ by following a suitable gradient flow, see for example Corollary 2.8 in \cite{mal1} (with minor adaptations). Somehow, the absence of critical points of $J_\rho$ above the level $b$ prevents us from having \emph{obstructions} while following the flow, see Figure \ref{fig:alti}.

\begin{pro} \label{contract}
Suppose $\rho=(\rho_1,\rho_2)\notin \L$. Then, there exists $b>0$ sufficiently large such that the sublevel $J_\rho^b$ is a deformation retract of $H^1(\Sg)$. In particular, it is contractible.
\end{pro}

The aim will be then to show how rich is the topological structure of the very low sublevels of $J_\rho$ and apply the Morse inequalities of Theorem \ref{morse} to deduce the main result of Theorem \ref{teo}. 
 
\

\subsection{Morse Theory}

\

\medskip

\no We recall here some classical results from Morse theory, which will be the main tool in proving Theorem \ref{teo}. 

Letting $N$ be a Hilbert manifold, we recall first that a function $f\in C^2(N,\R)$ is called a Morse function if all its critical points are nondegenerate. Moreover, the number of negative eigenvalues of the Hessian matrix at a critical point is called the index of the critical point. If $a<b$ are regular values of $f$ we then define the following sets:
\begin{equation} \label{set}
	\begin{array}{c}
		C_q(a,b) = \# \bigr\{  \mbox{critical points of $f$ in $\{a \leq f \leq b\}$ with index $q$}  \bigr\}, \vspace{0.3cm} \\ 
		\b_q(a,b) = \mbox{rank } \bigr( H_q \bigr( \{ f \leq b \}, \{ f \leq a \} \bigr) \bigr).
	\end{array}
\end{equation}
For the proof of the following result we refer for example to Theorem 4.3 in \cite{ch}.
\begin{thm}(\cite{ch}) \label{morse}
Let $N$ be a Hilbert manifold and $f\in C^2(N,\R)$ be a Morse function  satisfying the Palais-Smale condition. Let $a<b$ be regular values of $f$ and $C_q(a,b), \b_q(a,b)$ be as in \eqref{set}. Then the (respectively) strong and weak Morse inequalities hold true:
$$
	\sum_{q=0}^n (-1)^{n-q} C_q(a,b) \geq \sum_{q=0}^n (-1)^{n-q} \b_q(a,b), \qquad n=0,1,2,\dots
$$
$$
	C_q(a,b) \geq \b_q(a,b), \qquad n=0,1,2,\dots
$$
\end{thm}
The strategy will be to apply this result in our framework, namely with $N = H^1(\Sg)$ and $f=J_\rho$. We point out that the Palais-Smale condition is not necessarily needed for the Theorem \ref{morse} to hold, in fact it can be replaced by appropriate deformation lemmas for $f$, see Lemma 3.2 and Theorem 3.2 in \cite{ch}. The validity of such deformation lemmas can be obtained by following the ideas in \cite{mal1}, where a gradient flow for the scalar case \eqref{eq2} is defined.

For what concerns the assumption of $f$ to be a Morse function, one can repeat (with minor adaptations) the argument in \cite{demar}, which relies on a transversality result from \cite{sa-te}, to obtain the following result (recall the definition of $\M^2$ given in Subsection \ref{not}):
\begin{pro}(\cite{demar}) \label{non-deg}
Suppose $\rho=(\rho_1,\rho_2)\notin \L$. Then, for $(g,h)$ in an open dense subset of $\M^2 \times C^{\,2}(\Sg)^+$, $J_\rho$ is a Morse function.
\end{pro}
By the above discussion it follows that we are in position to apply Theorem \ref{morse} in our setting.

\section{Proof of the Main result}

\no We have now all the tools in order to prove the main result of Theorem \ref{teo}. Since the high sublevels of $J_\rho$ are contractible, see Proposition \ref{contract}, the goal will be to describe the topology of the low sublevels. 

\begin{figure}[h]
\centering
\includegraphics[width=0.3\linewidth]{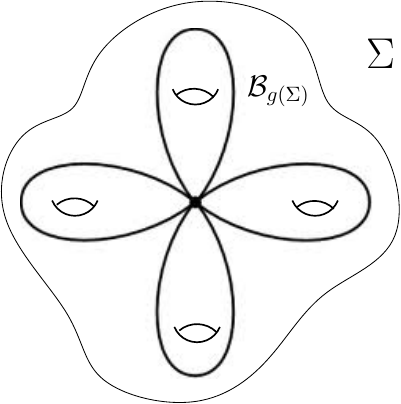}
\caption{}
\label{fig:bouquet}
\end{figure}

This will be done by means of a bouquet of circles and its homology will give then a bound of the number of solutions to \eqref{eq} by Theorem \ref{morse}.

\begin{figure}[h]
\centering
\includegraphics[width=0.7\linewidth]{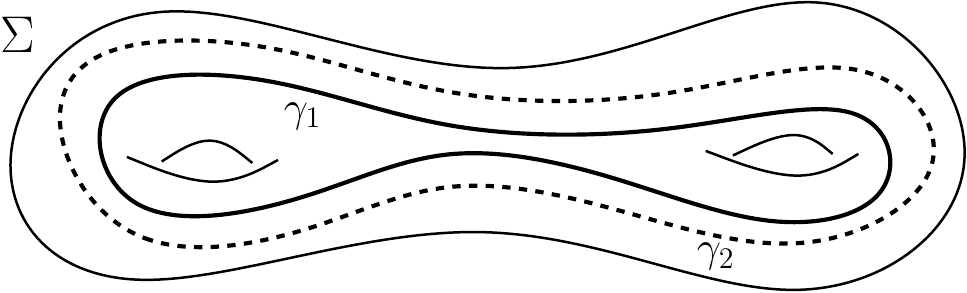}
\caption{}
\label{fig:curve}
\end{figure}

We recall that a bouquet $\mathcal B_N$ of $N$ circles (see Figure \ref{fig:bouquet}) is defined as $\mathcal B_N = \cup_{i=1}^N \mathcal S_i$, where $\mathcal S_i$ is homeomorphic to $S^1$ and $\mathcal S_i \cap \mathcal S_j = \{c\}$, and $c$ is called the center of the bouquet. The first simple result we need is the following, see the proof of Proposition 3.1 in \cite{bar-mar-mal}:
\begin{lem} \label{bouquet}
Let $\Sg$ be a surface with $g(\Sg)>0$. Then, there exist two curves $\g_1,\g_2 \subseteq \Sg$ satisfying (see Figure \ref{fig:curve})
\begin{itemize}
	\item[(1)] $\g_1$ and $\g_2$ do not intersect each other;
	\item[(2)] each of $\g_1$ and $\g_2$ are homeomorphic to respectively two disjoint bouquets of $g(\Sg)$ circles, see Figure \ref{fig:bouquet};
	\item[(3)] there exist global retractions $\Pi_i : \Sg \to \g_i$, $i=1,2$.
\end{itemize}
\end{lem}

We will now exploit the analysis developed in \cite{todatori} to describe the topology of the low sublevels of the functional $J_\rho$. As mentioned in the Introduction, by means of improved Moser-Trudinger inequalities one can deduce that if $\rho_1 < 8(k + 1)\pi$ and $\rho_2 < 8(l + 1)\pi$, then either $e^{u}$ is close to $\Sg_k$ or $e^{-u}$ is close to $\Sg_l$. This alternative is then expressed using the notion the topological join of $\Sg_k$ and $\Sg_l$, see \eqref{join}. Finally, applying the retractions $\Pi_1, \Pi_2$ introduced in the above Lemma, low energy sublevels may be described in terms of $(\g_1)_k * (\g_2)_l$ only. 

In fact, one can project the low sublevels of $J_\rho$ onto the latter set, see the proof of Proposition 4.7 and Section 6 in \cite{todatori}: for $\rho_1 \in (8k\pi, 8(k + 1)\pi)$, $\rho_2 \in (8l\pi, 8(l + 1)\pi)$ and for $L$ sufficiently large there exists a continuous map
$$
	\Psi : J_\rho^{-L} \to (\g_1)_k * (\g_2)_l.
$$
One the other hand, it is possible to do the converse, mapping $(\g_1)_k * (\g_2)_l$ into the low sublevels using suitable test functions, see Proposition 6.3 in \cite{todatori}:
$$
	\Phi : (\g_1)_k * (\g_2)_l \to J_\rho^{-L}.
$$
The above maps are somehow natural in the description of the low sublevels as we have the following important result, see Proposition 4.7 and Section 6 in \cite{todatori}:
\begin{thm}(\cite{todatori})
Suppose $\rho_1 \in (8k\pi, 8(k + 1)\pi)$, $\rho_2 \in (8l\pi, 8(l + 1)\pi)$ and $L$ sufficiently large. Then, the composition of the above maps $\Phi$ and $\Psi$ is homotopically equivalent to the identity map on $(\g_1)_k * (\g_2)_l$, i.e. $\Phi \circ \Psi \simeq \emph{Id}_{(\g_1)_k * (\g_2)_l}$.
\end{thm}
By the latter homotopy equivalence we directly deduce that the homology groups of $(\g_1)_k * (\g_2)_l$ are mapped injectively into the homology groups of $J_\rho^{-L}$ through the map induced by $\Phi$. 
\begin{cor}
Suppose $\rho_1 \in (8k\pi, 8(k + 1)\pi)$, $\rho_2 \in (8l\pi, 8(l + 1)\pi)$ and $L$ sufficiently large. Then, for any $q \in \N$ we have
$$
	H_q\bigr((\g_1)_k * (\g_2)_l\bigr) \hookrightarrow H_q\left(J_\rho^{-L}\right).
$$
\end{cor}
As a consequence of the above result we obtain a bound of the number of solutions to \eqref{eq} by Theorem \ref{morse}. One has just to observe that by Proposition \ref{contract}, taking $L\geq b$, the sublevel $J_\rho^L$ is contractible and therefore, by the long exact sequence of the relative homology, it follows that
$$
	H_{q+1}\left( J_\rho^L, J_\rho^{-L} \right) \cong \wtilde H_q \left(J_\rho^{-L}\right), \quad q\geq 0,
$$
$$	
	H_0\left( J_\rho^L, J_\rho^{-L} \right) = 0,
$$
where $\wtilde H_q(X)$ of a topological set $X$ is defined in the Subsection \ref{not}. Recalling the definition of $\b_q(a,b)$ introduced in \eqref{set} and the notation of $\wtilde \b_q$ given in the Subsection \ref{not}, the next result holds true by the above discussion and by taking $a=-L$ in Theorem \ref{morse}.
\begin{pro}\label{pro1}
Suppose $\rho_1 \in (8k\pi, 8(k + 1)\pi)$, $\rho_2 \in (8l\pi, 8(l + 1)\pi)$ and $L$ sufficiently large. Then, for any $q \in \N$ it holds that
$$
	\b_{q+1}(L,-L) \geq \wtilde\b_q\bigr((\g_1)_k * (\g_2)_l\bigr).
$$
\end{pro}
The next step is then to compute the homology groups of the topological join $(\g_1)_k * (\g_2)_l$. We recall that the two curves $\g_1$ and $\g_2$ were chosen such that there are homeomorphic to respectively two disjoint bouquets, see Lemma \ref{bouquet}. The homology group of the barycenters over this object was computed in Proposition 3.2 of \cite{bar-mar-mal}.
\begin{pro}(\cite{bar-mar-mal}) \label{pro2}
Let $\mathcal B_N$ be a bouquet of $N$ circles. Then, we have that
$$
	\wtilde H_q\bigr( (\mathcal B_N)_j \bigr) \cong \left\{ \begin{array}{lr}
																										     \Z^{\binom{j+N-1}{N-1}} & \mbox{if} \quad q=2N-1,  \\
																										     0                       & \mbox{if} \quad q\neq2N-1.
																									    \end{array} \right.
$$
\end{pro}
Finally, it is well known that the homology groups of the topological join of two sets $A$ and $B$ are expressed in terms of the sum of the homology groups of each set, see \cite{hat}.
\begin{pro}(\cite{hat}) \label{pro3}
Given two topological sets $A$ and $B$ we have
$$
	\wtilde H_q(A*B) \cong \bigoplus_{i=0}^q \wtilde H_i(A) \otimes \wtilde H_{q-i-1}(B).
$$
In particular it holds that
$$
 \wtilde\b_q(A*B) = \sum_{i=0}^q \wtilde\b_i(A) \, \wtilde\b_{q-i-1}(B).
$$
\end{pro}
We are now in position to deduce the main Theorem \ref{teo}. The proof will follow by applying the weak Morse inequality stated in Theorem \ref{morse} jointly with Proposition~\ref{pro1} and Propositions \ref{pro2}, \ref{pro3}. More precisely we get
$$
	\#\bigr\{ \mbox{solutions of \eqref{eq}}\bigr\}   \geq C_{q+1}(L,-L) \stackrel{Thm \, \ref{morse}}{\geq}  \b_{q+1}(L,-L) \stackrel{Prop \, \ref{pro1}}{\geq} \wtilde\b_q\bigr((\g_1)_k * (\g_2)_l\bigr)  \vspace{0.2cm}
$$
$$	
	          \hspace{-0.8cm}      \stackrel{Prop \, \ref{pro2} \, + \, Prop \, \ref{pro3}}{\geq}   \binom{k+g(\Sg)-1}{g(\Sg)-1} \binom{l+g(\Sg)-1}{g(\Sg)-1}   
$$
and the proof is concluded.

\

\begin{center}
\textbf{Acknowledgements}
\end{center}

\no Gratitude is expressed to  Professor Andrea Malchiodi for his support and for his kind help in preparing this paper.

The author is supported by the PRIN  project \emph{Variational and perturbative aspects of nonlinear differential problems}.

\end{document}